\documentclass{amsart}

\usepackage{lipsum}
\usepackage{amsfonts,bm,amssymb,amsmath}
\usepackage{graphicx}
\usepackage{epstopdf}
\usepackage[caption=false]{subfig}
\usepackage{pgfplots}
\usepackage{algorithmic}
\usepackage{amsopn}
\usepackage{amsrefs}
\usepackage{float}
\usepackage{geometry}

\usepackage{hyperref}
\usepackage{cleveref}

\usepackage{color}
\usepackage{diagbox}

\newtheorem{thm}{Theorem}[section]

\hypersetup{
	colorlinks,
	linkcolor={red!50!black},
	citecolor={blue!50!black},
	urlcolor={blue!80!black}
}

\numberwithin{equation}{section}

\definecolor{newcolor1}{rgb}{.8,.349,.1}
\colorlet{bblue}{blue!50!black}
\crefformat{equation}{(#2#1#3)}
\crefmultiformat{equation}{(#2#1#3)}{ and~(#2#1#3)}{, (#2#1#3)}{ and~(#2#1#3)}

\crefformat{figure}{Figure~#2#1#3}
\crefmultiformat{figure}{Figures~ #2#1#3}{ and~#2#1#3}{, (#2#1#3)}{ and~(#2#1#3)}
\crefformat{table}{Table~#2#1#3}

\def\f{\mbox{\boldmath $f$}}

\def\m{\mbox{\boldmath $m$}}

\def\p{\mbox{\boldmath $p$}}

\def\x{\mbox{\boldmath $x$}}

\def\0{\mbox{\boldmath $0$}}

\begin{document}

\title[A third order accurate scheme for LLG Equation]{Convergence analysis of a third order semi-implicit projection method for Landau-Lifshitz-Gilbert equation}





\author{Changjian Xie}
\address{School of Mathematics and Physics\\ Xi'an Jiaotong-Liverpool University\\ Ren'ai Rd. 111, Suzhou, 215123, Jiangsu\\ China.}
\email{Changjian.Xie@xjtlu.edu.cn}

\author{Cheng Wang}
\address{Mathematics Department\\ University of Massachusetts\\ North Dartmouth\\ MA 02747\\ USA.}
\email{cwang1@umassd.edu}


\date{\today}



\begin{abstract}	
The convergence analysis of a third-order scheme for the highly nonlinear Landau-Lifshitz-Gilbert equation with a non-convex constraint is considered. In this paper, we first present a fully discrete semi-implicit method  for solving the Landau-Lifshitz-Gilbert equation based on the third-order backward differentiation formula and the one-sided extrapolation (using previous time-step numerical values). A projection step is further used to preserve the length of the magnetization. We provide a rigorous convergence analysis for the fully discrete numerical solution by the introduction of two sets of approximated solutions where one set of solutions solves the Landau-Lifshitz-Gilbert equation and the other is projected onto the unit sphere. Third-order accuracy in time and fourth order accuracy in space is obtained provided that the spatial step-size is the same order as the temporal step-size and slightly large damping parameter $\alpha$ (greater than $\sqrt{2}/2$). And also, the unique solvability of the numerical solution without any assumption for the step-size in both time and space is theoretically justified, using a monotonicity analysis. All these theoretical results are guaranteed by numerical examples in both 1D and 3D spaces.
\end{abstract}

\maketitle

\section{Introduction}

The Landau-Lifshitz-Gilbert (LLG) equation is given by 
\begin{align} \label{equation-LL-alt}
{\m}_t=-{\m}\times\Delta{\m} + \alpha \Delta \m + \alpha | \nabla \m |^2 \m . 
\end{align}
with
\begin{equation}\label{boundary}
\frac{\partial{\m}}{\partial \boldmath {\nu}}\Big|_{\Gamma}=0,
\end{equation}
where $\Gamma = \partial \Omega$ and $\boldmath {\nu}$ is the unit outward normal vector along $\Gamma$. Here ${\m}\,:\,\Omega\subset\mathbb{R}^d\to S^2$ represents the magnetization vector field with $|{\m}|=1,\;\forall x\in\Omega$, $d=1,2,3$ is the spatial dimension, and $\alpha>0$ is the damping parameter. The first term on the right hand side of \cref{equation-LL-alt} is the gyromagnetic term, and the second and third terms are the damping term. Compared to the original LL equation \cite{Landau2015On}, \cref{equation-LL-alt} only includes the exchange term which poses the main difficulty in numerical analysis, as done in the literature \cite{weinan2001numerical, cimrak2004iterative, bartels2006convergence, gao2014optimal}. 
To ease the presentation, we set $\Omega=[0, 1]^d$, in which $d$ is the dimension. 


We have proposed a second order accurate finite difference scheme for model \eqref{equation-LL-alt} in previous work \cite{Xie20201JCP,Chen2021APNUM},  where we adopted the second-order Backward Differentiation Formula (BDF2) and the one-sided extrapolation (using previous time-step numerical values). A projection step is used to preserve the length of the magnetization. In this work, we propose and analyze a third-order accurate scheme that satisfies these desired properties. The third-order BDF approximation is applied to obtain an intermediate magnetization $\tilde{\m}$, and
the right-hand-side nonlinear terms are treated in a semi-implicit style with a
third-order extrapolation applied to the explicit coefficients and the lower order nonlinear terms. Such a numerical
algorithm leads to a linear system of equations with variable coefficients to solve at each time step. We provide a
rigorous convergence and error estimate, by the usage of the linearized stability analysis for the numerical error functions. In particular, we notice that, an \textit{a priori} $W_h^{1,\infty}$ bound assumption for the numerical solution at the previous time steps has h
to be imposed to pass through the convergence analysis. As a result, the
standard $L^2$ error estimate is insufficient to recover such a bound for the numerical solution. Instead, we have to perform the $H^1$ error estimate, and such a $W_h^{1,\infty}$ bound could be obtained at the next time step as a consequence of the $H^1$ estimate, via the help of the inverse inequality combined with a mild time step-size condition $k=\mathcal{O}(h)$ and the damping parameter constraint $\alpha>\sqrt{2}/2$. Careful error estimates for both the original magnetization $\m$ and the intermediate magnetization $\tilde{\m}$ have to be taken into consideration at the projection step (a highly nonlinear operation).

The rest of this paper is organized as follows. In \Cref{sec:main theory}, we introduce the fully discrete numerical scheme and state the main theoretical results. The proof of optimal rate convergence analysis in details are also provided. Numerical results are presented in \Cref{sec:experiments}, including both the 1-D and 3-D examples to confirm the theoretical analysis. Conclusions are drawn in \Cref{sec:conclusions}.

\section{Main theoretical results}
\label{sec:main theory}

\subsection{Numerical scheme}\label{discretisations}

In more details, the following numerical scheme is proposed: 
\begin{align}
\hat{\m}_h^{n+3} &= 3 \m_h^{n+2} - 3 \m_h^{n+1} + \m_h^n , \, \, \, 
\hat{\tilde{\m}}_h^{n+3} = 3 \tilde{\m}_h^{n+2} - 3 \tilde{\m}_h^{n+1} + \tilde{\m}_h^n , 
\label{cc} 
\end{align}

\begin{align}\label{scheme-1-1}
  & 
\frac{\frac{11}{6} \tilde{\m}_h^{n+3} - 3 \tilde{\m}_h^{n+2} + \frac32 \tilde{\m}_h^{n+1} 
- \frac13 \tilde{\m}_h^n }{k} \\ 
= &  - \hat{\m}_h^{n+3} \times \Delta_{h, (4)} \tilde{\m}_h^{n+3} 
 + \alpha \Delta_{h, (4)} \tilde{\m}_h^{n+3}  
 + \alpha | \tilde{\nabla}_{h, (4)} \hat{\m}_h^{n+3} |^2 \hat{\m}_h^{n+3} , \nonumber
\end{align}
\begin{align}
\m_h^{n+3} &= \frac{\tilde{\m}_h^{n+3}}{ |\tilde{\m}_h^{n+3}| } , \label{scheme-1-2}
\end{align}
%
%
%

\subsection{Main theoretical results}

The first theoretical result is the unique solvability analysis of scheme~\cref{cc}-\cref{scheme-1-2}.  We observe that the unique solvability for \cref{scheme-1-1} could be simplified as the analysis for
\begin{equation}\label{scheme-alt-1}
\frac{\frac{11}{6} \tilde{\m}_h - \p_h}{k}
=  - \hat{\m}_h \times \Delta_{h, (4)} \tilde{\m}_h
+ \alpha \Delta_{h, (4)}  \tilde{\m}_h + \f_h  , 
\end{equation}
with $\p_h$, $\hat{\m}_h$, $\f_h$ given.

\begin{thm} \label{thm:solvability}
	Given $\p_h$, $\hat{\m}_h$, $\f_h$, the numerical scheme \cref{scheme-alt-1} is uniquely solvable.
\end{thm}

The second theoretical result is the optimal rate convergence analysis. 
(Such theoretical result can be proved in a similar manner of our work \cite{xie2025}. The difference is that error analysis of the gyromagnetic term and the constraint of damping parameter. The details for the proof have been completed. However, the paper will be polished further.)

\begin{thm}\label{cccthm2} Let $\m_e \in C^4 ([0,T]; C^0) \cap C^3([0,T]; C^1) \cap L^{\infty}([0,T]; C^6)$ be the exact solution of \cref{equation-LL-alt} with the initial data $\m_e ({\x},0)=\m_e ^0({\x})$ and ${\m}_h$ be the numerical solution of the equation~\cref{cc}-\cref{scheme-1-2} with the initial data ${\m}_h^0=\m_{e,_h}^0$, $\m_h^1= \m _{e,h}^1$ and $\m_h^2= \m _{e,h}^2$. Suppose that the initial error satisfies $\|\m_{e,h}^\ell - \m_h^\ell \|_2 +\|\nabla_h ( \m_{e,h}^\ell - \m_h^\ell ) \|_2 = \mathcal{O} (k^3 + h^4),\,\ell=0,1, 2$, and $k\leq \mathcal{C}h$. In addition, we assume that $\alpha > \frac{\sqrt{2}}{2}$. Then the following convergence result holds as $h$ and $k$ goes to zero:
	\begin{align} \label{convergence-0} 
	\| \m_{e,h}^n - \m_h^n \|_{2}+ \|\nabla_h ( \m_{e,h}^n - \m_h^n ) \|_{2} &\leq \mathcal{C}(k^3+h^4) , \quad \forall n \ge 3 ,
	\end{align}	
	in which the constant $\mathcal{C}>0$ is independent of $k$ and $h$.
\end{thm}

\section{Numerical examples}
\label{sec:experiments}

In this section, we verify its accuracy in one-dimentional (1D) and three-dimentional (3D) cases.
In 1D, we choose the exact solution as below,
\begin{align*}
\m_e=[\cos(\cos(\pi x)) \sin (t), \sin(\cos(\pi x)) \sin(t),\cos(t)].
\end{align*}

For our proposed semi-implicit scheme in 1D, the spatial accuracy is shown in \cref{spaceAccuracy-v3}. The temporal accuracy is shown in \cref{timeAccuracy-v3}. 

\begin{table}[htbp]
	\centering
	\caption{Spatial accuracy for our scheme in 1D with $\alpha=0.01$, $N_t=1$e5, and $T=0.1$.}\label{spaceAccuracy-v3}
	\begin{tabular}{c|c|c|c}
		\hline
		$h$&$\|\m_h-\m_e\|_{\infty}$ &$\|\m_h-\m_e\|_2$ &$\|\m_h-\m_e\|_{H^1}$ \\
		\hline
		$1/16$&9.081724873419295e-06 & 6.480755332622150e-06 & 9.925805477628938e-05 \\
		\hline
		$1/32$& 5.803494078082672e-07 & 4.113288938373731e-07 & 6.397213459888633e-06  \\
		\hline
		$1/64$&3.641971156598256e-08 & 2.580904089084249e-08 & 4.030617118174344e-07\\
		\hline
		$1/128$&  2.280198946325029e-09 & 1.614584426499681e-09 & 2.524313256978843e-08 \\
		\hline
		$1/256$&  1.424844547903703e-10 & 1.009176797442278e-10 & 1.578646012187375e-09\\
		\hline
		order &3.99		 &3.99&3.99\\
		\hline
	\end{tabular}
\end{table}

\begin{table}[htbp]
	\centering
	\caption{Temporal accuracy for our scheme in 1D with $\alpha=0.01$, $N_x=1$e4, and $T=0.1$.}\label{timeAccuracy-v3}
	\begin{tabular}{c|c|c|c}
		\hline
		$k$&$\|\m_h-\m_e\|_{\infty}$ &$\|\m_h-\m_e\|_2$ &$\|\m_h-\m_e\|_{H^1}$ \\
		\hline
		$T/8$ & 7.159874484963247e-09 & 4.515028250584236e-09 & 1.600888891412794e-08 \\
		\hline
		$T/12$& 2.271219642913103e-09 & 1.419508208447264e-09 & 5.068087073023508e-09\\
		\hline
		$T/16$& 9.832080288818545e-10 & 6.118815615963404e-10 & 2.191042705781703e-09 \\
		\hline
		$T/24$& 2.946392574365575e-10 & 1.836536334121961e-10 & 6.673805158610665e-10\\
		\hline
		$T/32$ & 1.260517107359860e-10 & 7.822989610437513e-11 & 2.801951682042707e-10\\
		\hline
		order &2.92		 &2.93&2.92\\
		\hline
	\end{tabular}
\end{table}

In 3D, we take the exact solution as below,
	\begin{align*}
\m_e=[\cos(\cos(\pi x)\cos(\pi y)\cos(\pi z)) \sin (t), \sin(\cos(\pi x)\cos(\pi y)\cos(\pi z)) \sin(t),\cos(t)].
\end{align*}
The results for spatial accuracy are presented in \Cref{spaceAccuracy-D-3}. The temporal accuracy is shown in \Cref{timeAccuracy-D-1}.

\begin{table}[htbp]
	\centering
	\caption{Spatial accuracy for our scheme in 3D with $\alpha=0.01$, $N_t=1$e4 and $T=1$.}\label{spaceAccuracy-D-3}
	\begin{tabular}{c|c|c|c}
		\hline
		$h$&$\|\m_h-\m_e\|_{\infty}$ &$\|\m_h-\m_e\|_2$ &$\|\m_h-\m_e\|_{H^1}$ \\
		\hline
		$1/4$& 0.139274466278363& 0.065240602569682  &0.184903434480020\\ 
		\hline
		$1/6$& 0.032104973270598 & 0.014122998077254 & 0.042420430109726  \\
		\hline
		$1/8$& 0.010931127155401 & 0.004685104278569 & 0.014199407041049  \\
		\hline
		$1/10$& 0.004653055687641 & 0.001968267343486 & 0.005962002221913 \\
		\hline
		$1/12$ &0.002294249564032 & 9.633572776318690e-04 & 0.002911324475901\\
		\hline
		order &	3.74	 & 3.84&  3.78\\
		\hline
	\end{tabular}
\end{table}

\begin{table}[htbp]
	\centering
	\caption{Temporal accuracy for our scheme in 3D with $\alpha=0.01$, and $T=1$.}\label{timeAccuracy-D-1}
	\begin{tabular}{c|c|c|c|c}
		\hline
		$h$&$k,\;k^3\approx h^4$&$\|\m_h-\m_e\|_{\infty}$ &$\|\m_h-\m_e\|_2$ &$\|\m_h-\m_e\|_{H^1}$ \\
		\hline
		$1/6$&$1/10$& 0.013506716901339 & 0.006419417702565 & 0.028945176405124\\
		\hline
		$1/8$&$1/15$& 0.004867787436378 & 0.002150405678805 & 0.010031179319299 \\
		\hline
		$1/10$&$1/21$& 0.001527458902012 & 6.625474139793397e-04 & 0.003637833068190 \\
		\hline
		$1/12$&$1/27$&6.898974024560633e-04 & 3.004612076000572e-04 & 0.001706824873478  \\
		\hline
		$1/16$&$1/40$&2.064831665434252e-04 & 8.654322306801727e-05 & 5.302069678519973e-04  \\
		\hline
		order&-- &3.06 &3.14&2.90\\
		\hline
	\end{tabular}
\end{table}

\section{Conclusions}
\label{sec:conclusions}

In this paper, we have proposed and analyzed a third-order time stepping scheme to solve the LLG equation. The third-order BDF is applied for temporal discretization and a linearized multistep approximation is used for the nonlinear terms. The resulting scheme avoids a well-known difficulty associated with the nonlinearity of the system, and its unique solvability is established via the monotonicity analysis of the system. In addition, an optimal rate convergence analysis is provided, by making use of a linearized stability analysis for the numerical error functions, in which the $W_h^{1,\infty}$ error estimate at the projection step has played an important role. Numerical experiments in both 1D and 3D cases are presented to verify the unconditional stability and the third-order accuracy in time and fourth order accuracy in space.

\section*{Acknowledgments}

This work is supported in part by the grants NSF DMS-2012669 (C.~Wang), and Jiangsu Science
and Technology Programme-Fundamental Research Plan Fund BK20250468, Research and Development Fund 
of XJTLU (RDF-24-01-015) (C. Xie).

\bibliographystyle{amsplain}
\bibliography{draft1_BDF3}
%
%
%
%

\end{document}